\documentclass{amsart}

\usepackage{epsfig}
\usepackage{graphicx}
\usepackage{color}
\usepackage{amsthm,amsmath,amsfonts,amssymb,amscd}
\usepackage{lineno,hyperref,amssymb,amsmath}

\usepackage{enumerate}
\def\Circ{{\sf A}}

\def\Circ{{\rm Circ}}\def\AA{{\sf A}}

\def\GG{{\sf {G}}}

\def\JJ{{\sf {J}}}

\def\Circ{{\sf {M}}}
\def\PP{{\sf {P}}}

\def\Circ{\hbox{{\scriptsize{\sl a}}}}

\def\1{{\sf {1}}}
\def\0{{\sf {0}}}
\def\a{{\sf {a}}}
\def\b{{\sf {b}}}
\def\c{{\sf {c}}}

\def\e{{\sf {e}}}

\def\g{{\sf {g}}}
\def\h{{\sf {h}}}

\def\t{{\sf {t}}}
\def\u{{\sf {u}}}
\def\w{{\sf {w}}}
\def\v{{\sf {v}}}

\def\qed{\hfill {\hbox{\footnotesize{$\Box$}}}}
\def\0{{\sf {\sf 0}}}

\def\a{{\sf a}}
\def\k{{\sf k}}

\def\w{{\sf w}}

\def\v{{\sf v}}
\def\0{{\sf 0}}

\def\so{{\sf supp}}

\def\proof{{\noindent\bf Proof.}\hskip 0.3truecm}
\def\BBox{\kern  -0.2cm\hbox{\vrule width 0.15cm height 0.3cm}}

\def\Circ{{\rm Circ}}
\def\AA{{\sf A}}

\def\B{\mathcal{B}}
\def\E{\mathcal{E}}
\def\C{\mathcal{C}}

\def\F{\mathcal{F}}
\def\G{\mathcal{G}}

\def\K{\mathcal{K}}
\def\L{\mathcal{L}}
\def\M{\mathcal{M}}
\def\P{\mathcal{P}}

\def\GG{{\sf G}}
\def\II{{\sf I}}
\def\JJ{{\sf J}}

\def\PP{{\sf P}}

\def\NN{\mathbb{N}}
\def\ZZ{\mathbb{Z}}

\def\RR{\mathbb{R}}

\def\demo{{\noindent\sc Proof.}\hskip 0.3truecm}
\def\BBox{\kern  -0.2cm\hbox{\vrule width 0.15cm height 0.3cm}}

\oddsidemargin .cm \evensidemargin 1.cm \textwidth=16.5cm
\include {mak}
\begin{document}

% Your \newcommans below (if there are any):
\newtheorem{propo}{Proposition}[section]
\newtheorem{lemma}[propo]{Lemma}
\newtheorem{theorem}[propo]{Theorem}
\newtheorem{corollary}[propo]{Corollary}
\newtheorem{definition}[propo]{Definition}

\def\1{{\sf 1}}
\def\0{{\sf 0}}
\def\e{{\sf e}}
\def\a{{\sf a}}
\def\v{{\sf v}}
\def\y{{\sf y}}
\def\z{{\sf z}}
\def\g{{\sf g}}
\def\h{{\sf h}}
\def\II{{\sf I}}
\def\JJ{{\sf J}}
\def\so{{\sf supp}}
\def\proof{{\noindent  Proof.}\hskip 0.3truecm}
\def\BBox{\kern  -0.2cm\hbox{\vrule width 0.15cm height 0.3cm}}
\def\L{{\mathcal  L}}\def\P{{\mathcal  P}}
\def\B{{\mathcal B}}
\def\C{{\mathcal C}}
\def\K{{\mathcal K}}\def\M{{\mathcal M}}
\def\E{{\mathcal E}}
\def\F{{\mathcal F}}
\def\G{{\mathcal G}}
\oddsidemargin 16.5mm
\evensidemargin 16.5mm

\thispagestyle{plain}

\hspace{-.5cm}
\includegraphics[scale=0.5]{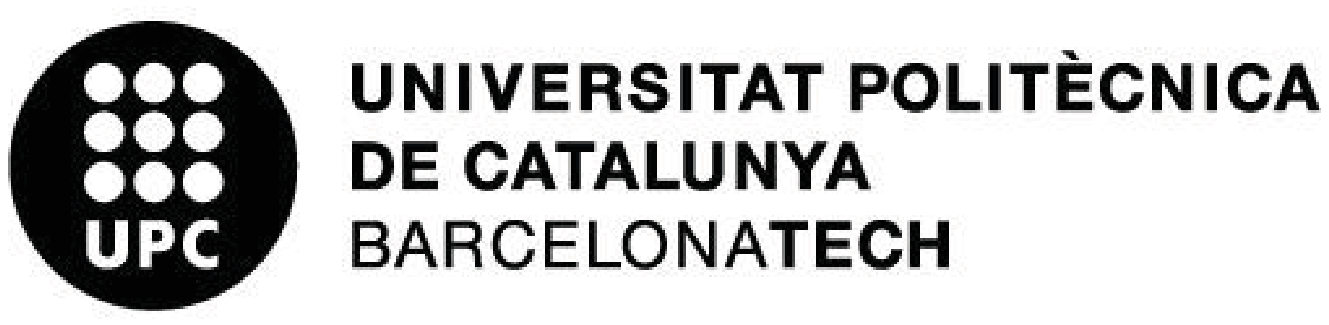} \hspace{3cm}
\includegraphics[scale=0.14]{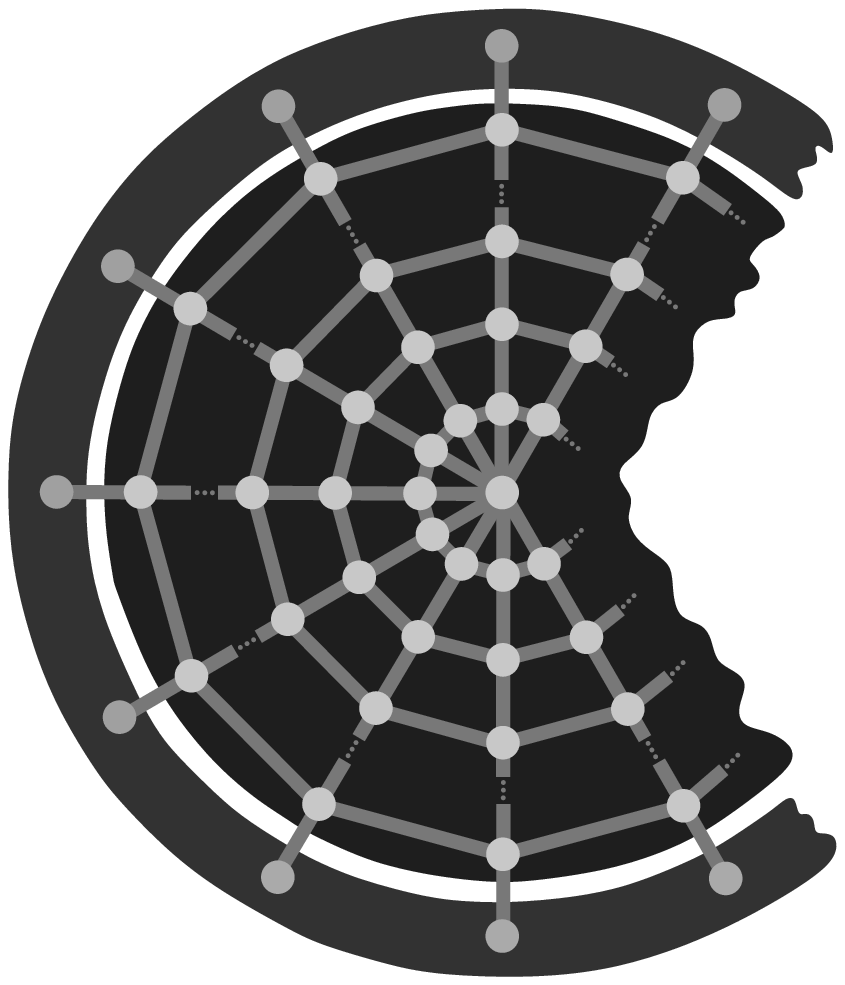}
\vspace{-1.3cm}

\hfill\begin{minipage}{.28\textwidth}
\hspace{.47cm}\noindent {\small\bf GRUPO MAPTHE}

%\noindent {\small\sc Dept. MAIII}

\hspace{.5cm}\noindent {\scriptsize \bf Preprint P2-2015}

\hspace{.5cm}\noindent {\scriptsize \bf \today.}
\end{minipage}
%\noindent{\scriptsize Available electronically at {\bf http://www-ma3.upc.es/users/bencar}}

\vspace{5cc}
\begin{center}

{\Large\bf  The inverse matrix of some circulant matrices
%INVERSAS DE ALGUNAS MATRICES\\ SIM\'ETRICAS Y CIRCULANTES
\rule{0mm}{6mm}\renewcommand{\thefootnote}{}%Enter at least one, but not more than 4 MSCs.
% First entered MSC will be a primary one, others (at most 3) will be secondary.
\footnotetext{\scriptsize 2000 Mathematics Subject Classification:  }}

\vspace{1cc} {\large\it A. Carmona, A.M. Encinas, S. Gago, M.J. Jim\'{e}nez, M. Mitjana }

\vspace{1cc}
\parbox{24cc}{{\small{\bf Abstract.}
We present here  necessary and sufficient conditions for the invertibility of  circulant and symmetric  matrices that depend on three parameters and moreover, we   explicitly compute the inverse. The  techniques we use are related with the solution of boundary value problems associated to second order linear difference equations. Consequently, we reduce the computational cost of the problem. In particular, we recover the inverses of some well known circulant matrices whose coefficients are  arithmetic or geometric sequences,  Horadam numbers among others. We also characterize when a general symmetric circulant and tridiagonal matrix is invertible and in this case, we compute explicitly its inverse.
}}
\end{center}

\vspace{1cc}

%\vspace{1.5cc}

\section{Introduction}
\label{introduction}
Many problems in applied mathematics and science lead to the solution of linear systems having circulant coefficients related to the periodicity of the problems, as the ones that appear  when using the finite difference method to approximate elliptic equations with periodic boundary conditions, see \cite{ChCh92}. Circulant matrices have a wide range of application in signal processing, image processing, digital image disposal, linear forecast, error correcting code theory, see  \cite{CoVaCl04, WaDo07}. %\cite{CVC04, WD07}. 
 In the last years, there have been several papers on circulant matrices that attend to give an effective expression for the determinant, the eigenvalues and the inverse of the matrix, see for instance \cite{Fu11,ShCeHa11, YaTa13}.%\cite{F11,SCH11, YT13}.

In this paper, we consider circulant matrices of type $\Circ(a,b,c,\ldots,c)$ and $\Circ(a,b,c,\ldots,c,b)$. This type of matrices raise when dealing, for example, with finite differences for solving one dimensional elliptic equations, or when computing the Green function of some networks obtained by the addition of new vertices to a previously known one, see \cite{CaEnGaMi15}. We give a necessary and sufficient condition for its invertibility. Moreover, as it is known, their inverse is a circulant matrix and we explicitly obtain a closed formula for the expression of the coefficients.

For fixed $n\in \NN^*$, we consider the vector space $\mathbb{R}^n$ together with the standard inner product $\langle \cdot,\cdot\rangle$. Moreover, we  denote the components of the vector $\v\in \RR^n$ as $v_j$, $j=1,\ldots,n$, \textit{i.e.}, $\v=(v_1,\ldots,v_n)^\top$.

As usual,  $\1$ is the all ones vector and  $\0$ is the all zeroes  vector  in $\RR^n$. In addition, $\e$ is the $\RR^n$ vector whose first component is equal to 1 and 0 otherwise.

The set of  matrices of order $n$ with real coefficients is denoted $\M_n(\RR)$.  Moreover,
matrix $\II$ is  the Identity matrix of order $n$ and $\JJ$ is the matrix of order $n$ whose entries are  all ones.

A matrix $\AA=(a_{ij})$  is named {\it circulant  with  parameters $a_1,\ldots,a_n$} if
\begin{equation}
\label{circ}
\AA=\begin{bmatrix}a_1 & a_2 & \cdots & a_n\\
a_n & a_1 & \cdots & a_{n-1}\\
\vdots & \vdots & \ddots & \vdots\\
a_{2} & a_3 & \cdots & a_1\end{bmatrix}\end{equation}
or equivalently,
\begin{equation}
\label{circ:coef}
%a_{ij}=\left\{\begin{array}{cl}a_{j+1-i}, & i\le j,\\ a_{n+1+j-i}, & i>j,\end{array}\right.
a_{ij}=a_{1+(j-i)({\rm mod\,}n)}\end{equation}
see \cite{ShCeHa11, YaTa13}.%\cite{SCH11, YT13}.

Given $\a\in \RR^n$,  $\Circ(\a)=\Circ(a_1,\ldots,a_n)\in \M_n(\RR)$ is the  {\it circulant matrix} with parameters $a_1,\ldots,a_n$.
Notice that $\Circ(\e)=\II$ and  $\Circ(\1)=\JJ$.

Let $\tau$ be a  permutation  of the set $\{1,\ldots,n\}$  defined as,
\begin{equation}
\label{tau}
\tau(1)=1,\hspace{.25cm}\tau(j)=n+2-j,\hspace{.25cm}j=2,\ldots,n.
\end{equation}

We denote $\PP_\tau\in \M_n(\RR)$ the   matrix with entries $(p_{ij})$, such that for any  $j=1,\ldots,n$, $p_{\tau(j)j}=1$ and $p_{ij}=0$, otherwise. It turns out that  $\PP_\tau$ is invertible  and  it is satisfied $\PP_\tau^{-1}=\PP^\top_\tau=\PP_{\tau^{-1}}=\PP_\tau$.

In a similar way we now define  $\a_\tau=\PP_\tau\a$; that is,   the vector whose components are  $(a_\tau)_1=a_1$ and $(a_\tau)_j=a_{n+2-j}$, $j=2,\ldots,n$. Thus,  $\1_\tau=\1$ and $\langle \a_\tau,\1\rangle=\langle \a,\1\rangle.$ Moreover,
\begin{equation}
\label{tau:ident}
\Circ(\a_\tau)=\PP_\tau \Circ(\a)\PP_\tau.
\end{equation}

For any $\a\in \RR^n$,  the matrices
\begin{eqnarray}
\label{circ:sigma}
\Circ_\tau(\a)=\PP_\tau\Circ(\a)=&  \begin{bmatrix}a_1 & a_2 & \cdots & a_n\\
a_2 & a_3 & \cdots & a_{1} \\
\vdots & \vdots & \ddots & \vdots\\
a_{n} & a_1 & \cdots & a_{n-1}\end{bmatrix}\hspace{.25cm}\hbox{and}\\
 \Circ^\tau(\a)=\Circ(\a)\,\PP_\tau=& \hspace{-.75cm}\begin{bmatrix}a_1 & a_{n} & \cdots & a_2\\
a_n & a_{n-1} & \cdots & a_{1}\\
\vdots & \vdots & \ddots & \vdots\\
a_{2} & a_1 & \cdots & a_3\end{bmatrix}\end{eqnarray}
are named {\it left--circulant and right--circulant with parameters $a_1,\ldots,a_n$},  respectively. Both matrices are symmetric  and for this reason  matrices   $\Circ_\tau(\a)$ were called circulant symmetric matrices in \cite{WaDo07}. %\cite{WD07}.
To avoid confusion, we will refer to this type of matrices with the notation introduced above. In addition,  from identity \eqref{tau:ident} we get $\Circ^\tau(\a)=\Circ_\tau(\a_\tau)$ for any $\a\in \RR^n$.

For the sake of completeness, we enumerate the properties of the circulant matrices that are relevant in our study. All the statements are easily shown to hold.

\begin{lemma}
\label{circ:prop}
Given any  $\a\in \RR^n$, the following properties hold:
\begin{itemize}
\item[(i)]   For any $\b\in \RR^n$ and $\alpha,\beta\in \RR$, $\Circ(\alpha\a+\beta \b)=\alpha\Circ(\a)+\beta \Circ(\b)$.
\item[(ii)]   $\Circ(\a)^\top=\Circ(\a_{\tau}) $. In particular, $\Circ(\a)$ is  symmetric iff $ \a=\a_{\tau}.$
\item[(iii)] $\Circ_\tau(\a)=\Circ( \a)$ iff $\Circ^\tau(\a)=\Circ( \a)$. This equalities hold iff $\a=\a_\tau$.
\item[(iv)]   $\Circ(\a)\1=\langle \a,\1\rangle\,\1$. Moreover, if $\Circ(\a)$ is invertible  then $\langle \a,\1\rangle\not=0$.
\item[(v)]   For any $\b\in \RR^n$, $\Circ(\a)\b=\Circ(\b_{\tau}) \a_{\tau}$ and $\Circ(\a)\Circ(\b)=\Circ(\b)\Circ(\a)=\Circ(\c_{\tau})$, where $\c=\Circ(\a) \b_{\tau}=\Circ( \b) \a$.
\item[(vi)] $\Circ(\a)$ is invertible iff the linear system   $\Circ(\a)\g=\e$ is compatible. In that case, there exists a only solution  $\g(\a)$ that, in addition satisfies $\langle \g(\a),\1\rangle =\langle \a,\1\rangle^{-1}$. Moreover, $\Circ(\a)^{-1}=\Circ\big(\g(\a)\big)^\top$ and  $\a_\tau=\a$ iff $ \g(\a)_\tau=\g(\a)$.
\item[(vii)] $\Circ_\tau(\a)$ and $\Circ^\tau(\a)$ are invertible iff $\Circ(\a)$ is invertible and, in that case is
$$\Circ_\tau(\a)^{-1}=\Circ_\tau\big(\g(\a)\big) \hspace{.25cm}\hbox{and}\hspace{.25cm}\Circ^\tau( \a)^{-1}=\Circ^\tau\big(\g(\a)\big).$$
\end{itemize}
\end{lemma}
\vspace{.25cm}

One of the main problems in the field of circulant matrices is to determine invertibility conditions  and, in this case,   to calculate the inverse. From the previous Lemma,  it turns out that for any $\a \in\mathbb{R}^n$  the invertibility of matrices $\Circ_\tau(\a)$ and $\Circ^\tau(\a)$, as well as the computation of their inverses, can be deduced from  the invertibility and the inverse matrix of  $\Circ(\a)$.  Furthermore, the solution of this last problem can be reduced to study the compatibility of  a certain linear system  and the computation of its solution, when it exists. Moreover, the problem has been widely studied in the literature by using the primitive  $n$--th root of unity and some polynomial associated with it,  see \cite{Gr01, ShCeHa11}.
%\cite{G01, SCH11}.
Specifically, let  $\omega=e^{\frac{2\pi}{n}i}$ be the primitive  $n$--th root of unity. In addition,  define  for each $ j = 0, \ldots, n-1$,  the vector $ \t_j = \big(1,\omega^j, \ldots, \omega ^{j(n-1)} \big)^{\top} \in \RR ^ n $ and  for any $ \a \in \RR ^n $   the polynomial $P_\a(x)=\sum\limits_{j=1}^na_jx^{j-1}$. Observe that $ \t_0 = \1$ and for any $\a \in \RR^n $, $P_\a(1)=\langle \a,\1\rangle$. The following lemma provides a necessary and sufficient condition for the invertibility of $\Circ(\a)$ and gives a formula for the inverse.

\begin{lemma}
\label{circ:propinv}
For any $\a\in \RR^n$, the following properties hold:
\begin{itemize}
\item[(i)]  $\Circ(\a)\t_j=P_{\a}(\omega^j)\t_j$, for any $j=0,\ldots,n$. In particular, ${\rm det\,}\Circ(\a)=\prod\limits_{k=0}^{n-1}P_\a(\omega^k)$.
\item[(ii)]  $\Circ(\a)$ is invertible iff $P_\a(\omega^j)\not=0$, $j=0,\ldots,n-1$. In this case, $\Circ^{-1}(\a)=\Circ( \h_\a)$ where $({ h_\a})_{_j}=\dfrac{1}{n}\sum\limits_{k=0}^{n-1}\omega^{-k(j-1)}P_\a(\omega^k)^{-1}$.
\end{itemize}
\end{lemma}

Let us notice that property (i) of the previous Lemma implies that all circulant matrices of order $n$ have the same eigenvectors but different  eigenvalues. Although the  problem  is completely solved, the computational complexity of formula (ii) for the determination of the inverse of a circulant invertible matrix grows with the order of the matrix, so it is not useful from the computational point of view. To illustrate this phenomena, we apply Lemma \ref{circ:propinv} to the cases $n=2,3$, which can be directly solved  without  difficulty.

For $n=2$,  $\omega=-1$ and for a given  $\a\in \RR^2$, is $P_\a(x)=a_1+a_2x$, and $P_\a(-1)=a_1-a_2$. Thus,  $\Circ(\a)$ is invertible iff $(a_1+a_2)(a_1-a_2)\not=0$  and
$$
\begin{array}{rl}
({h_\a})_{_1}=&\hspace{-.25cm}\dfrac{1}{2}\left(\dfrac{1}{a_1+a_2}+\dfrac{1}{a_1-a_2}\right)=\dfrac{a_1}{a_1^2-a_2^2},\\[3ex]
({h_\a})_{_2}=&\hspace{-.25cm}\dfrac{1}{2}\left(\dfrac{1}{a_1+a_2}-\dfrac{1}{a_1-a_2}\right)=\dfrac{a_2}{a_1^2-a_2^2}.
\end{array}
$$
On the other hand, if $n=3$,  $\omega=\frac{1}{2}(-1+i\sqrt{3})$, $\omega^2=\overline \omega$  and for a given $\a\in \RR^3$, is $P_\a(x)=a_1+a_2x+a_3x^2$;   $P_\a(\omega)=a_1+a_2\omega+a_3\overline \omega$ and  $P_\a(\omega^2)=\overline{P_\a(\omega)}$. Thus,   $\Circ(\a)$ is invertible iff
$$P_\a(1)|P_\a(\omega)|^2=(a_1+a_2+a_3)|a_1+a_2\omega+a_3\overline \omega|^2=a_1^3+a_2^3+a_3^3-3a_1a_2a_3\not=0$$
in this case,
$$\begin{array}{rl}
({h_\a})_{_1}\!\!\!=&\hspace{-.25cm}\dfrac{1}{3P_\a(1)|P_\a(\omega)|^2}\Big(|P_\a(\omega)|^2+2P_\a(1)\Re(P_\a(\omega)\Big)\!\!=\!\!\dfrac{a_1^2-a_2a_3}{a_1^3+a_2^3+a_3^3-3a_1a_2a_3},
\\[3ex]
({h_\a})_{_2}\!\!\!=&\hspace{-.25cm}\dfrac{1}{3P_\a(1)|P_\a(\omega)|^2}\Big(|P_\a(\omega)|^2+2P_\a(1)\Re(\omega P_\a(\omega)\Big)\!\!=\!\!\dfrac{a_3^2-a_1a_2}{a_1^3+a_2^3+a_3^3-3a_1a_2a_3},
\\[3ex]
({h_\a})_{_3}\!\!\!=&\hspace{-.25cm}\dfrac{1}{3P_\a(1)|P_\a(\omega)|^2}\Big(|P_\a(\omega)|^2+2P_\a(1)\Re(\overline \omega P_\a(\omega)\Big)\!\!=\!\!\dfrac{a_2^2-a_1a_3}{a_1^3+a_2^3+a_3^3-3a_1a_2a_3}.
\end{array}$$
\vspace{.5cc}

%Nuestro prop\'osito es obtener  resultados expl\'icitos para la determinaci\'on de la inversas de algunas matrices circulantes sim\'etricas de order $n\ge 3$, que surgen habitualmente en las aplicaciones,\footnote{detallar algunas, como por ejemplo el c\'alculo de la resistencia efectiva y del \'indice de Kirchhoff de cadenas lineales.}
%reduciendo  por tanto el coste computacional de aplicar la f\'ormula del Lema anterior.

We aim here to compute the inverse matrix of some circulant matrices of order $n\ge 3$ with three parameters at most. For instance, this kind of circulant matrices  appears when computing the effective resistance and the Kirchhoff index of a network that comes from the addition of new nodes to a previously known one, see \cite{CaEnGaMi15}. We significantly reduce the computational cost of applying Lemma \ref{circ:propinv} since the key point for finding the mentioned inverse matrix consists in  solving a difference equation of order at most two.

\section{Matrices $\Circ(a,b,c,\ldots,c)$}

Given $a,b,c\in \RR$, let $\a(a,b,c)\in \RR^n$ defined as $\a(a,b,c)=(a,b,c,\ldots,c)^\top$.  Then,  $\Circ(a,b,c,\ldots,c)=\Circ\big(\a(a,b,c)\big)$.

For any $q\in \RR$, we also consider  the vector $\z(q)=(q^{n-1},q^{n-2},\ldots,q,1)^\top\in \RR^n$. It immediately follows that  $\z_\tau(q)=(q^{n-1},1,q,\ldots,q^{n-2})^\top$ and  $\langle \z(q),\1\rangle=\dfrac{q^n-1}{q-1}$. We observe that the last identity  also holds for   $q=1$, since $\lim\limits_{q\to 1}\dfrac{q^n-1}{q-1}=n=\langle \z(1),\1\rangle$. We remark that $\z(1)=\1$ and $\z_\tau(-1)=\z(-1)=(-1,1,-1,\ldots,1)^\top$, for $n$ even.

\begin{propo}
\label{first}
For any  $q\in \RR$ it is satisfied
$$\Circ\big(\a(q,-1,0)\big) \z_\tau(q)=[q^n-1]\e.$$
Moreover, the following properties hold:
\begin{itemize}
\item[(i)] $\Circ\big(\a(q,-1,0)\big)$ is invertible iff $q^n\not=1$,  and the inverse matrix is
 $$\Circ\big(\a(q,-1,0)\big)^{-1}=(q^n-1)^{-1}\Circ\big( \z(q)\big).$$
\item[(ii)] The linear system $\Circ\big(\a(1,-1,0)\big)\h=\v$ is compatible iff $\sum\limits_{i=1}^nv_i=0$. In this case, for any $\gamma\in \RR$ the unique solution of the linear system  satisfying   $\langle \h,\1\rangle=\gamma$ is
$$h_j=\dfrac{1}{n}\Big[\gamma-\sum\limits_{i=1}^niv_i\Big]+\sum\limits_{i=j}^{n}v_i,\hspace{.25cm}j=1,\ldots,n.$$
\item[(iii)] If $n$ is even, the linear system $\Circ\big(\a(-1,-1,0)\big)\h=\v$ is compatible iff $\langle \v,\z(-1)\rangle=0$. In this case, all the solutions are given by
$$h_j=(-1)^{j+1}\sum\limits_{i=j}^n(-1)^iv_i+\alpha\z(-1),\hspace{.25cm}j=1,\ldots,n, \hspace{.25cm}\alpha\in\mathbb{R}$$
and hence  $2\langle \h,\1\rangle=\langle \v,\1\rangle$.
 \end{itemize}
\end{propo}

The inverse matrix of a circulant matrix whose parameters are a geometric sequence can be computed as an application of the above result.

\begin{corollary}
\label{geometric}
For any $a,r\in \RR$, the matrix   $\Circ(ar^{n-1},\ldots,ar,a)$ is invertible iff $a(r^n-1)\not=0$. In this case, the inverse matrix is
$$\Circ(ar^{n-1},\ldots,ar,a)^{-1}=\big(a(r^n-1)\big)^{-1}\Circ\big(\a(r,-1,0)\big).$$
\end{corollary}
\vspace{.25cm}

The main result on the present section displays a necessary and sufficient condition for the invertibility of the matrices we are dealing with. Moreover, when the inverse exists, we can provide a simple and closed expression for its entries.

\begin{theorem}
\label{three parameters:first}
For any $a,b,c\in \RR$,  the circulant matrix $\Circ(a,b,c,\ldots,c)$ is invertible iff
$$\big[a+b+(n-2)c\big]\big[(a-b)^2+\big(1-(-1)^n\big)(c-b)^2\big]\not=0$$
and, in that case,
$\Circ(a,b,c,\ldots,c)^{-1}=\Circ\big(\k(a,b,c)\big)$ where, if $a\not=2c-b$

$$k_j(a,b,c)=\dfrac{(c-b)^{j-1}(a-c)^{n-j}}{\big((a-c)^n-(c-b)^n)}-\dfrac{c}{(a+b-2c)\big(a+b+(n-2)c\big)},$$
j=1,\ldots,n, and
$$k_j(2c-b,b,c)=\dfrac{1}{n}\Big[\frac{1}{nc}+\dfrac{n-1}{2(c-b)}\Big]-\dfrac{(j-1)}{n(c-b)},\hspace{.25cm}j=1,\ldots,n.$$
\end{theorem}
\demo Define $\alpha(n;a,b,c)=\big[a+b+(n-2)c\big]\big[(a-b)^2+\big(1-(-1)^n\big)(c-b)^2\big]$, thus $\alpha(n;a,b,b)=0$ iff $a=b$ or $a=(1-n)b$. Moreover,  when $c\not=b$, if  $n$ is odd,  then $\alpha(n;a,b,c)=0$ iff $a=-b-(n-2)c$, while if    $n$ is even, then $\alpha(n;a,b,c)=0$ iff either $a=b$ or  $a=-b-(n-2)c$.

According to  (iv) in Lemma \ref{circ:prop}, $\langle \a(a,b,c),\1\rangle=a+b+(n-2)c\not=0$ is a necessary condition in order to ensure the invertibility of $\Circ\big(\a(a,b,c)\big)$. So, from now on we will assume that this condition holds.
Moreover, from Lemma \ref{circ:prop} (vi),  $\Circ\big(\a(a,b,c)\big)$ is invertible iff the linear system   $\Circ\big(\a(a,b,c)\big)\h=\e$ is    compatible, and in that case there is a unique solution that, in addition, satisfies   $\langle \h,\1\rangle=\langle \a(a,b,c),\1\rangle^{-1}$.

Therefore, since $\Circ\big(\a(a,b,c)\big)=\Circ\big(\a(a-c,b-c,0)\big)+c\JJ$, if $\h\in \RR^n$ fulfills  $\Circ\big(\a(a,b,c)\big)\h=\e$, then
$\Circ\big(\a(a-c,b-c,0)\big)\h=\e-c\langle \a(a,b,c),\1\rangle^{-1}\1$.

Reciprocally, if $\h\in \RR^n$ satisfies  $\Circ\big(\a(a-c,b-c,0)\big)\h=\e-c\langle \a(a,b,c),\1\rangle^{-1}\1$,  then
$$\Circ\big(\a(a,b,c)\big)\h=\Circ\big(\a(a-c,b-c,0)\big)\h+c\JJ\h=\e+c\big[\langle \h,\1\rangle-\langle \a(a,b,c),\1\rangle^{-1}\big]\,\1$$
and hence,   $\h$ is a solution of $\Circ\big(\a(a,b,c)\big)\h=\e$ iff $\langle \h,\1\rangle=\langle \a(a,b,c),\1\rangle^{-1}$.

Consequently, we have shown that  
$$\Circ\big(\a(a,b,c)\big)\h=\e\hspace{.15cm}\hbox{iff}\hspace{.15cm}\Circ\big(\a(a-c,b-c,0)\big)\h=\e-c\langle \a(a,b,c),\1\rangle^{-1}\1$$
and moreover $\langle \h,\1\rangle=\langle \a(a,b,c),\1\rangle^{-1}$.

If $c=b$, then   $\langle \a(a,b,b),\1\rangle=a+(n-1)b$ and $\Circ\big(\a(a-c,b-c,0)\big)=(a-b)\II$. Thus, the system $(a-b)\II\h=\e-b\big(a+(n-1)b\big)^{-1}\1$ is compatible iff $a\not=b$ and then
$$\h=\dfrac{1}{(a-b)\big(a+(n-1)b\big)}\Big[\big(a+(n-1)b\big)\e-b\1\Big],$$
which implies $\langle \h,\1\rangle=\dfrac{\big(a+(n-1)b\big)-bn}{(a-b)\big(a+(n-1)b\big)}=\dfrac{1}{a+(n-1)b}=\langle \a(a,b,b),\1\rangle^{-1}$. Notice that,
$$
\begin{array}{rl}
h_1(a,b,b)=& \hspace{-.25cm}\dfrac{1}{a-b}-\dfrac{b}{(a-b)\big(a+(n-1)b\big)}, \\[3ex]
h_j(a,b,b)=& \hspace{-.25cm}-\dfrac{b}{(a-b)\big(a+(n-1)b\big)},\hspace{.25cm}j=2,\ldots,n.
\end{array}
$$
For the case $c\not=b$, we consider $q=\dfrac{a-c}{c-b}$, then $\a(a-c,b-c,0)=(c-b)\a(q,-1,0)$ and, as a consequence,  the system $\Circ\big(\a(a-c,b-c,0)\big)\h=\e-c\langle \a(a,b,c),\1\rangle^{-1}\1$ is equivalent to
$$\Circ\big(\a(q,-1,0)\big)\h=\dfrac{1}{(c-b)\big(a+b+(n-2)c\big)}\Big(\big(a+b+(n-2)c\big)\e-c\1\Big).$$

If $\h$ is a solution of the previous system,  then
$$\begin{array}{rl}
\dfrac{(a+b-2c)}{(c-b)\big(a+b+(n-2)c\big)}=&\hspace{-.25cm}\langle \Circ\big(\a(q,-1,0)\big)\h,\1\rangle=\langle \h,\Circ\big( \a_\tau(q,-1,0)\big)\1\rangle\\[1ex]
 =&\hspace{-.25cm}\langle \a(q,-1,0),\1\rangle\langle \h,\1\rangle=\dfrac{(a+b-2c)}{(c-b)}\,\langle \h,\1\rangle,\end{array}$$
which implies that if $a+b-2c\not=0,$  then  $\langle \h,\1\rangle\!=\!\dfrac{1}{a+b+(n-2)c}\!=\!\langle \a(a,b,c),\1\rangle^{-1}$. So that,  if $c\not=b$ and $a+b-2c\not=0$; i.e. $q\not=1$, then
$
\Circ\big(\a(a,b,c)\big)\h=\e$
 iff \hspace{.5cm}$\Circ\big(\a(q,-1,0)\big)\h=\dfrac{1}{(c-b)\big(a+b+(n-2)c\big)}\Big(\big(a+b+(n-2)c\big)\e-c\1\Big).$

For $n$ odd or $n$ even but $b\not=a$; i.e. $q^n\not=1$, according to  Proposition \ref{first} (i), $\Circ\big(\a(q,-1,0)\big)$ is invertible and in addition,
 $$\begin{array}{rl}\h=&\hspace{-.25cm}\dfrac{1}{(c-b)(q^n-1)\big(a+b+(n-2)c\big)}\Circ\big(\z(q)\big)\Big(\big(a+b+(n-2)c\big)\e-c\1\Big)\\[3ex]
 =&\hspace{-.25cm}\dfrac{1}{(c-b)(q^n-1)\big(a+b+(n-2)c\big)}\Big(\big(a+b+(n-2)c\big) \z_\tau(q)-c\langle\z(q),\1\rangle\,\1\Big).
 \end{array}$$

If  $n$ is even and $b=a$, then $q=-1$, $a+b+(n-2)c=2(b-c)+nc$  and
$$\langle \big(2b+(n-2)c\big)\e-c\1, \z(-1)\rangle=-\big(2(b-c)+nc\big)\not=0.$$
 Consequently, the system $$\Circ\big(\a(-1,-1,0)\big)\h=\dfrac{1}{(c-b)\big(a+b+(n-2)c\big)}\Big(\big(a+b+(n-2)c\big)\e-c\1\Big)$$ is incompatible, thus $\Circ\big(\a(a,b,c)\big)$ has no inverse.

If  $c\not=b$ and $a+b-2c=0$, i.e. $q=1$, then $a+b+(n-2)c=nc$ and the system
$$\Circ\big(\a(1,-1,0)\big)\h=\dfrac{1}{n(c-b)}(n\e-\1)$$
has solution. Moreover, by Proposition  \ref{first} (ii), if
$$h_1=\dfrac{1}{n}\Big[\frac{1}{nc}+\dfrac{n-1}{2(c-b)}\Big],\hspace{.5cm}h_j=\dfrac{1}{n}\Big[\frac{1}{nc}+\dfrac{n-1}{2(c-b)}\Big]-\dfrac{(n+1-j)}{n(c-b)},\hspace{.25cm}j=2,\ldots,n$$
then,  $\h$ is the only solution of the linear system that satisfies $\langle \h,\1\rangle=\dfrac{1}{nc}$.

In all cases  it suffices to take $\k=\h_\tau$.\qed

 Notice that  for $a=2c-b$ in the previous theorem, $\Circ\big(\a(2c-b,b,c\big)$ is invertible iff $c(c-b)\not=0$ and then,
 the entries of vector $\h(2c-b,b,c)$ are the elements of an arithmetic sequence. Thus we can give a characterization of the inverse matrix of a circulant matrix whose parameters are in arithmetic sequence. Of course, our results coincide with those obtained in \cite{BaSo10}
 %\cite{BS10}
 and, for left--circulant matrices with parameters in arithmetic progression, in \cite{WaDo07}.

\begin{corollary}
\label{arit}
For any $a,b\in \RR$, the   matrix $\Circ\big(a,a+b,\ldots,a+(n-1)b\big)$ is invertible iff $\big(2a+(n-1)b\big)b\not=0$ and in this case,
$$\Circ\big(a,a+b,\ldots,a+(n-1)b\big)^{-1}=\dfrac{2}{n^2\big(2a+(n-1)b\big)}\,\JJ-\dfrac{1}{nb}\,\Circ\big(\a(1,-1,0)\big).$$
%i
In particular, for any $m\in\ZZ$ such that $2m+n\not=1$, the matrix $\Circ(m,m+1,\ldots,m+n-1)$ is invertible and its inverse is
$$\Circ(m,m+1,\ldots,m+n-1)^{-1}=\frac{2}{n^2(2m+n-1)}\,\JJ-\frac{1}{n}\Circ\big(\a(1,-1,0)\big).$$
\end{corollary}

\section{Matrices $\Circ(a,b,c,\ldots,c,b)$}

For any $a,b,c\in \RR$, let  $\b(a,b,c)\in \RR^n$ defined as $\b(a,b,c)=(a,b,c,\ldots,c,b)$. Then, $\Circ(a,b,c,\ldots,c,b)=\Circ\big(\b(a,b,c)\big)$ and $ \b_\tau(a,b,c)=\b(a,b,c)$, since matrix $\Circ(a,b,c,\ldots,c,b)$ is symmetric. Regarding the case  $\b(a,b,b)=\a(a,b,b)$,  matrix $\Circ(a,b,b,\ldots,b,b)$ has been analyzed in the previous section, so from now on we assume  $c\not=b$.
The case $c=0$ has been analyzed in \cite{Ro90}
% \cite{R90}
 under the name of symmetric circulant tridiagonal matrix, assuming the condition $|a|>2|b|>0$; that is, that $\Circ\big(\b(a,b,0)\big)$ is a strictly diagonally dominant matrix.

 Notice that  $\Circ\big(\b(2,-1,0)\big)$ is nothing but the so called {\it combinatorial Laplacian} of a $n$--cycle. More generally, for any $q\in \mathbb{R}$, $\Circ\big(\b(2q,-1,0)\big)$ is the matrix associated with the    Schr\"odinger operator on the cycle with constant potential $2(q-1)$ and hence its inverse is the Green's function of a $n$--cycle; or equivalenty, it can be seen as the Green function associated with a path with periodic boundary conditions, see \cite{BeCaEn09}.
%\cite{BCE09}.
Since the inversion of matrices of type $\Circ\big(\b(2q,-1,0)\big)$ involves  the resolution of second order difference equations with constant coefficients, we enumerate some of their properties.

 A \textit{Chebyshev sequence} is a sequence of polynomials $\{Q_n(x)\}_{n\in \ZZ}$ that satisfies the recurrence
\begin{equation}
\label{chebyshev}
Q_{n+1}(x)=2xQ_n(x)-Q_{n-1}(x),\hspace{.25cm}\hbox{for
each}\hspace{.15cm}n\in \ZZ.\end{equation}
Recurrence
(\ref{chebyshev}) shows that any Chebyshev sequence is uniquely
determined by the choice of the corresponding zero and one order
 polynomials, $Q_0$ and $Q_1$ respectively. In particular,  the sequences
\begin{math}\{T_n\}_{n=-\infty}^{+\infty}\end{math} and \begin{math}\{U_n\}_{n=-\infty}^{+\infty}\end{math}
 denote the \emph{first and second  kind
Chebyshev polynomials} that are obtained when we choose
\begin{math}T_0(x)=U_0(x)=1,\end{math} \begin{math}T_1(x)=x\end{math}, \begin{math}U_1(x)=2x\end{math}.

Next we describe some properties of the Chebyshev polynomials of first and second kind that will be useful in the present work. See \cite{MaHa03} %\cite{MH03}
for proofs and more details.

 \begin{itemize}
 \item[(i)] For any Chebyshev sequence
   $\{Q_n\}_{n=-\infty}^{+\infty}$ there exists $\alpha,\beta\in \RR$ such that $Q_n(x)=\alpha U_{n-1}(x)+\beta U_{n-2}(x)$, for any $n\in \ZZ$.
 \item[(ii)]  $T_{-n}(x)=T_n(x)$ and $U_{-n}(x)=-U_{n-2}(x)$, for any $n\in \ZZ$. In particular, $U_{-1}(x)=0$.
 \item[(iii)] $T_{2n+1}(0)=U_{2n+1}(0)=0$, $T_{2n}(0)=U_{2n}(0)=(-1)^n$, for any $n\in \ZZ$.
 \item[(iv)] Given  $n\in \NN^*$ then, $T_n(q)=1$ iff $q=\cos\left(\frac{2\pi j}{n}\right)$, $j=0,\ldots,\lceil\frac{n-1}{2}\rceil$, whereas $U_n(q)=0$ iff $q=\cos \big(\frac{\pi j}{n+1}\big)$, $j=1,\ldots,n$. In this case, $U_{n-1}(q)=(-1)^{j+1}$ and $U_{n+1}(q)=(-1)^j$.
% $$U_{n-1}(q)=\dfrac{\sin\big(n\frac{\pi j}{n+1}\big)}{\sin\big(\frac{\pi j}{n+1}\big)}=\dfrac{\sin\big(j\pi -\frac{\pi j}{n+1}\big)}{\sin\big(\frac{\pi j}{n+1}\big)}=-\cos(j\pi)=(-1)^{j+1}$$ y
% $$U_{n+1}(q)=\dfrac{\sin\big((n+2)\frac{\pi j}{n+1}\big)}{\sin\big(\frac{\pi j}{n+1}\big)}=\dfrac{\sin\big(j\pi+\frac{\pi j}{n+1}\big)}{\sin\big(\frac{\pi j}{n+1}\big)}=\cos(j\pi)=(-1)^j$$
 %
 \item[(v)]  $T_n(1)=1$ and  $U_n(1)=n+1$, whereas $T_n(-1)=(-1)^n$ and $U_n(-1)=(-1)^n(n+1)$, for any  $n\in \ZZ$. %
\item[(vi)] $T_n(x)=xU_{n-1}(x)-U_{n-2}(x)$ y $T_n'(x)=nU_{n-1}(x)$, for any $n\in \ZZ$.
\item[(vii)] $2(x-1)\sum\limits_{j=0}^nU_j(x)=U_{n+1}(x)-U_n(x)-1$, for any $n\in \NN$.
 \end{itemize}

Chebyshev recurrence \eqref{chebyshev}  encompasses all linear second order recurrences with constant coefficients, see \cite{AhBeDr05},
%\cite{ABD05},
 so we can consider more general recurrences. Let  $\{H_n(r,s)\}_{n=0}^\infty$, where $r,s\in \ZZ$ and $s\not=0$, the  {\it Horadam numbers} defined as the solution of the recurrence
\begin{equation}
\label{horadam}
H_{n+2}=rH_{n+1}+sH_n,\hspace{.25cm}H_0=0,\hspace{.15cm}H_1=1.
\end{equation}

Notice that for any $n\in \NN^*$, $H_n(1,1)=F_n$, the $n$--th {\it Fibonacci number}, $H_n(2,1)=P_n$, the $n$--th {\it Pell number}, $H_n(1,2)=J_n$, the $n$--th {\it Jacobsthal number} and $H_n(2,-1)=U_{n-1}(1)=n$.

The equivalence between any second order difference equation and Chebyshev equations leads to  the following result, see \cite[Theorem 3.1]{AhBeDr05} and \cite[Theorem 2.4]{EJ15}.
%\cite[Theorem 3.1]{ABD05} and \cite{EJ15}.
\begin{lemma}
\label{equiv}
Given $r,s\in \ZZ$ and $s\not=0$, we have the following results:
\begin{itemize}
\item[(i)] If $s<0$, then $H_n(r,s)=(\sqrt{-s})^{n-1}U_{n-1}\big(\frac{r}{2\sqrt{-s}}\big)$, $n\in \NN^*$.
\item[(ii)] If $s>0$, then $H_{2n}(r,s)=rs^{n-1}U_{n-1}\big(1+\frac{r^2}{2s}\big)$, $n\in \NN^*$.
\end{itemize}
In particular, for any $n\in \NN^*$, $F_{2n}=U_{n-1}\big(\frac{3}{2}\big)$, $J_{2n}=2^{n-1}U_{n-1}\big(\frac{5}{4}\big)$, and $P_{2n}=2U_{n-1}(3)$. In addition, $H_{2n}(r,r)=r^{n}U_{n-1}\big(1+\frac{r}{2 }\big)$ when $r>0$ and $H_n(r,r)=(\sqrt{-r})^{n-1}U_{n-1}\big(\frac{\sqrt{-r}}{2}\big)$ for $r<0$.
\end{lemma}
\vspace{.5cc}

 In addition, for any $q\in \RR$  we  denote by  $\u(q)$, $\v(q)$ and $\w(q)$ the vectors in $\RR^n$ whose components are   $u_j=U_{j-2}(q)$, $v_j=U_{j-1}(q)$ and $w_j=U_{j-2}(q)+U_{n-j}(q)$, respectively.

\begin{lemma}
\label{vector} For any $q\in \RR^n$, the following properties hold:
\begin{itemize}
\item[(i)] $ \w_\tau(q)=\w(q)$ and $\langle \w(q),\1\rangle=\dfrac{T_n(q)-1}{q-1}$. Moreover, $\w(1)=n\1$.
\item[(ii)] $\w(q)=\0$ iff $q=\cos\left(\frac{2\pi j}{n}\right)$, $j=1,\ldots,\lceil\frac{n-1}{2}\rceil$. In this case, $\langle \u(q),\1\rangle=\langle \v(q),\1\rangle=0$.
\item[(iii)] When $n$ is even, then $w_{2j-1}(0)=0$ and $w_{2j}(0)=(-1)^{j-1}\big[1-(-1)^{\frac{n}{2}}\big]$, $j=1,\ldots, \frac{n}{2}$.
\item[(iv)] When $n$ is odd, then $w_{2j-1}(0)=(-1)^{\frac{n+1}{2}+j}$, $j=1,\ldots, \frac{n+1}{2}$ and $w_{2j}(0)=(-1)^{j-1}$, $j=1,\ldots, \frac{n-1}{2}$.
\item[(iv)] When $n$ is odd, then $w_j(-1)=(-1)^{j-1}(n+2-2j)$, $j=1,\ldots,n$.
\end{itemize}
\end{lemma}
\demo $\w(q)=\0$ iff $U_{n-j}(q)=-U_{j-2}(q)$ for any $j=1,\ldots,n$ and this equality holds iff $U_{n-1}(q)=0$ and $U_{n-2}(q)=-1$.
Moreover,   $U_{n-1}(q)=0$ iff  $q=\cos\left(\frac{k\pi }{n}\right)$, $k=1,\ldots,n-1$, thus   $U_{n-2}(q)=(-1)^{k+1}$, leads to  $U_{n-2}(q)=-1$ iff $k=2j$. \qed

{\sc Remark:}  The quotient $\frac{T_n(q)-1}{q-1}$ is well defined for  $q=1$, because    $T_n(1)=1$, $U_n(1)=n+1$, and $T_n'(q)=nU_{n-1}(q)$, using  l'H\^{o}pital's rule, $\lim\limits_{q\to 1}\langle \w(q),\1\rangle=nU_{n-1}(1)=n^2$. Moreover,  for $q=1$, is $\w(1)=n\1$ thus,  $\langle \w(1),\1\rangle=n^2$.
\vspace{.5cc}

\begin{propo}
\label{cycle}
For any  $q\in \RR$,
$$\Circ\big(\b(2q,-1,0)\big)\w(q)=2[T_n(q)-1]\e.$$
and the following  holds:
\begin{itemize}
\item[(i)] $\Circ\big(\b(2q,-1,0)\big)$ is invertible iff $q\not=\cos\left(\frac{2\pi j}{n}\right)$, $j=0,\ldots,\lceil\frac{n-1}{2}\rceil$ and,
 $$\Circ\big(\b(2q,-1,0)\big)^{-1}=\frac{1}{2[T_n(q)-1]}\Circ\big(\w(q)\big).$$
\item[(ii)] If $q=1$, the linear system $\Circ\big(\b(2q,-1,0)\big)\h=\v$ is compatible iff $\langle \v,\1\rangle=0$ in this case, for any  $\gamma\in \RR$  the only  solution satisfying  $\langle \h,\1\rangle=\gamma$ is given by
$$h_j=\dfrac{\gamma}{n}-\dfrac{1}{2n}\sum\limits_{i=1}^n|j-i|(n-|i-j|)v_i,\hspace{.25cm}j=1,\ldots,n.$$
\item[(iii)]If $q=\cos\left(\frac{2\pi j}{n}\right)$, $j=1,\ldots,\lceil\frac{n-1}{2}\rceil$,  the linear system $\Circ\big(\b(2q,-1,0)\big)\h=\v$ is compatible iff $\langle \h,\u(q)\rangle=\langle \h,\v(q)\rangle=0$.
 \end{itemize}
\end{propo}
\demo To prove (i), notice that $\w(q)$ is the first column of the Green function for the  Schr\"odinger operator for a $n$--cycle, or equivalently for  a $(n+1)$--path with periodic boundary conditions, see \cite[Proposition 3.12]{BeCaEn09}.

To prove (ii), it suffices to see that $\GG=(g_{ij})$, where $g_{ij}=\dfrac{1}{12n}\big(n^2-1-6|i-j|(n-|i-j|)\big)$, $i,j=1,\ldots,n$ is the  Green  function of the Combinatorial Laplacian of the cycle, see for instance \cite{BeCaEnMi102}.
% \cite{BCEM10}.
 The third claim (iii), comes from   (ii) of Lemma \ref{vector} that states $\w(q)=\0$.  In addition,  in this case, $U_{n-1}(q)=0$, $U_{n-2}(q)=-1$ and $U_{n}(q)=1$.  Besides, vectors $\u(q)$ and $\w(q)$ satisfy
$$\begin{array}{rl}
2qu_1-u_2-u_n=& \hspace{-.25cm}-1-U_{n-2}(q)=0, \\[1ex]
 -u_1-u_{n-1}+2qu_n=& \hspace{-.25cm}-U_{n-3}(q)+2qU_{n-2}(q)=U_{n-1}(q)=0,\\[2ex]
2qv_1-v_2-v_n=& \hspace{-.25cm}2q-2q-U_{n-1}(q)=0,\\[1ex]
-v_1-v_{n-1}+2qv_n=& \hspace{-.25cm}-1-U_{n-2}(q)+2qU_{n-1}(q)=0,
\end{array}$$
thus,  $\Circ\big(\b(2q,-1,0)\big)\u(q)=\Circ\big(\b(2q,-1,0)\big)\v(q)=\0$.    \qed
\vspace{.25cm}

Next, the main result in this section is proved.  We give necessary and sufficient conditions for the existence of the inverse of matrix $\Circ(a,b,c,\ldots,c,b)$ and we explicitly obtain the coefficients of the inverse, when it exists.

\begin{theorem}
For $a,b,c\in \RR$, the circulant matrix $\Circ(a,b,c,\ldots,c,b)$ is invertible iff
$$\big(a+2b+(n-3)c\big)\prod\limits_{j=1}^{\lceil\frac{n-1}{2}\rceil}\Big[a-c+2(b-c)\cos\Big(\frac{2\pi j}{n}\Big)\Big]\not=0$$
and, in this case
$$\Circ(a,b,c,\ldots,c,b)^{-1}=\Circ\big(\g(a,b,c)\big),$$
where if $a\ne 3c-2b $
$$g_j(a,b,c)=\dfrac{U_{j-2}(q)+U_{n-j}(q)}{2(c-b)[T_n(q)-1]}-\dfrac{c}{(a+2b-3c)\big(a+2b+(n-3)c\big)},\hspace{.25cm}j=1,\ldots,n,$$
with
$q=\dfrac{c-a}{2(b-c)}$  , whereas
$$g_j(3c-2b,b,c)=\dfrac{1}{12n(c-b)}\big(n^2-1-6(j-1)(n+1-j)\big)+\frac{1}{n^2c},\hspace{.25cm}j=1,\ldots,n.$$
\end{theorem}
\demo From claim (iv) of  Lemma \ref{circ:prop},  a necessary condition for the invertibility of $\Circ\big(\b(a,b,c)\big)$ is $\langle \b(a,b,c),\1\rangle=a+2b+(n-3)c\not=0$, so, we will assume that this  condition holds.
Moreover, claim (vi) in the same Lemma states that a necessary and sufficient condition to  get $\Circ\big(\b(a,b,c)\big)$  invertible is   the compatibility of the linear system  $\Circ\big(\b(a,b,c)\big)\g=\e$, and in that case there is an only solution that satisfies  $\langle \g,\1\rangle=\langle \b(a,b,c),\1\rangle^{-1}$.

 As in  Theorem \ref{three parameters:first},   
 $$\Circ\big(\b(a,b,c)\big)\g=\e\hspace{.15cm}\hbox{iff}\hspace{.15cm}\Circ\big(\b(a-c,b-c,0)\big)\g=\e-c\langle \b(a,b,c),\1\rangle^{-1}\1$$ 
 and  moreover, $\langle \g,\1\rangle=\langle \b(a,b,c),\1\rangle^{-1}$.

Since $\b(a-c,b-c,0)=(c-b)\b(2q,-1,0)$, the linear system  
$$\Circ\big(\b(a-c,b-c,0)\big)\g=\e-c\langle \b(a,b,c),\1\rangle^{-1}\1$$ 
is equivalent to system
$$\Circ\big(\b(2q,-1,0)\big)\g=\dfrac{1}{(c-b)\big(a+2b+(n-3)c\big)}\Big(\big(a+2b+(n-3)c\big)\e-c\1\Big).$$

If $\g$ is a solution of the above system, then
$$\begin{array}{rl}
\dfrac{(a+2b-3c)}{(c-b)\big(a+2b+(n-3)c\big)}=&\hspace{-.25cm}\langle \Circ\big(\b(2q,-1,0)\big)\g,\1\rangle=\langle \g,\Circ\big(\b(2q,-1,0)\big)\1\rangle \\[1ex]
=&\hspace{-.25cm}\langle \b(2q,-1,0),\1\rangle\langle \g,\1\rangle=\dfrac{(a+2b-3c)}{(c-b)}\,\langle \g,\1\rangle.\end{array}$$
As a consequence, if $a+2b-3c\not=0$ then $\langle \g,\1\rangle=\dfrac{1}{a+2b+(n-3)c}=\langle \b(a,b,c),\1\rangle^{-1}$. Under this assumption; that is,   if $a\not=3c-2b$ or equivalently  $q\not=1$, then
$\Circ\big(\b(a,b,c)\big)\g=\e\hspace{.25cm}$ iff $$\Circ\big(\b(2q,-1,0)\big)\g=\dfrac{1}{(c-b)\big(a+2b+(n-3)c\big)}\Big(\big(a+2b+(n-3)c\big)\e-c\1\Big).$$

In addition, if $\prod\limits_{j=1}^{\lceil\frac{n-1}{2}\rceil}\Big[a-c+2(b-c)\cos\Big(\frac{2\pi j}{n}\Big)\Big]\not=0$, then $q\not=\cos\Big(\frac{2\pi j}{n}\Big)$, for any $j=1,\ldots,\lceil\frac{n-1}{2}\rceil$. Using claim  (i) in Proposition \ref{cycle}, $\Circ\big(\b(2q,-1,0)\big)$ is invertible, and
{\small{ $$\begin{array}{rl}\g=&\hspace{-.25cm}\dfrac{1}{2(c-b)\big(a+2b+(n-3)c\big)[T_n(q)-1]}\Circ\big(\w(q)\big)\Big(\big(a+2b+(n-3)c\big)\e-c\1\Big)\\[3ex]
 =&\hspace{-.25cm}\dfrac{1}{2(c-b)\big(a+2b+(n-3)c\big)[T_n(q)-1]}\Big(\big(a+2b+(n-3)c\big)\w(q)-c\langle\w(q),\1\rangle\,\1\Big).
 \end{array}$$}}

If  there exists $j=1,\ldots,\lceil\frac{n-1}{2}\rceil$, such that  $a-c+2(b-c)\cos\Big(\frac{2\pi j}{n}\Big)=0$, i.e.  $q=\cos\Big(\frac{2\pi j}{n}\Big)$, then, statement (ii) in Lemma \ref{vector} ensures
$$\langle \big(a+2b+(n-3)c\big)\e-c\1,\v(q)\rangle=\big(a+2b+(n-3)c\big) v_1(q)=a+2b+(n-3)c\not=0$$
so, by claim  (iii) in Proposition \ref{cycle}, the linear system
$\Circ\big(\b(a,b,c)\big)\g=\e$ is incompatible and,
$\Circ\big(\b(a,b,c)\big)$ is not invertible.

When $a=3c-2b$, this is $q=1$, then $a+2b+(n-3)c=nc$ and system
$$\Circ\big(\b(2,-1,0)\big)\g=\dfrac{1}{n(c-b)}(n\e-\1)$$
is  compatible. Moreover, using claim (ii) in Proposition \ref{cycle},  the vector   $\g\in \RR^n$ whose components  are given for any $j=1,\ldots,n$ by
$$g_j=\dfrac{1}{n^2c}-\dfrac{1}{2n(c-b)}(j-1)\big(n-(j-1)\big)+\dfrac{1}{2n^2(c-b)}\sum\limits_{i=1}^n|j-i|(n-|i-j|),$$
is the only solution  of the system satisfying    $\langle \g,\1\rangle=\dfrac{1}{nc}$. Last, we only have to  take into account that
$\sum\limits_{i=1}^n|j-i|(n-|i-j|)=\dfrac{n}{6}(n^2-1)$, for any $j=1,\ldots,n$. \qed
\vspace{.25cm}

 The case  $a=3c-2b$ in the above theorem, involves the Green function of a cycle. Cases related to this,  raise as application in the analysis of problems associated with this combinatorial structures.

\begin{corollary}
For a given $a,b\in \RR$,   matrix
$$\AA=\Circ\big(a,a+b(n-1),a+2b(n-2),\ldots,a+jb(n-j),\ldots,a+b(n-1)\big)$$
is invertible iff $\big(6a+b(n^2-1)\big)b\not=0$ and,
$$\AA^{-1}=\dfrac{6}{n^2\big(6a+b(n^2-1)\big)}\,\JJ-\dfrac{1}{2nb}\,\Circ\big(\b(2,-1,0)\big).$$
\end{corollary}

\begin{corollary}
For a given $a,b\in \RR$, the following results hold:
\begin{itemize}
\item[(i)] If $n=1\, {\rm mod}(4)$, then $\AA=\Circ(a,a,b,b,a,a,\ldots,a,a,b,b,a)$ is invertible iff $(a-b)\big(a(n+1)+b(n-1)\big)\not=0$ and then
\begin{center}
$\AA^{-1}=\dfrac{1}{a-b}\Circ(\b(0,1,0)\big)-\dfrac{2(a+b)}{(a-b)\big(a(n+1)+b(n-1)\big)}\,\JJ$
\end{center}
\item[(ii)] If $n=2\, {\rm mod}(4)$, then  $\AA=\Circ\big(\frac{a+b}{2},a,\frac{a+b}{2},b,\frac{a+b}{2},\ldots,\frac{a+b}{2},b,\frac{a+b}{2},a)$ is invertible iff $(a-b)\big(a(n+1)+b(n-1)\big)\not=0$ and then
\begin{center}
$\AA^{-1}=\dfrac{1}{a-b}\Circ(\b(0,1,0)\big)-\dfrac{2(a+b)}{(a-b)\big(a(n+1)+b(n-1)\big)}\,\JJ$\end{center}
\item[(iii)] If $n=3\, {\rm mod}(4)$, then $\AA=\Circ(b,a,a,b,b,\ldots,a,a,b,b,a,a)$ is invertible iff $(a-b)\big(a(n+1)+b(n-1)\big)\not=0$ and then
\begin{center}
$\AA^{-1}=\dfrac{1}{a-b}\Circ(\b(0,1,0)\big)-\dfrac{2(a+b)}{(a-b)\big(a(n+1)+b(n-1)\big)}\,\JJ$
\end{center}
\item[(iv)] When $n$ is odd, then 
$$\AA=\Circ\big(a+nb,a-(n-2)b,\ldots,a+(-1)^{j-1}(n+2-2j)b,\ldots,a-(n-2)b\big)$$ 
is invertible iff $b(an+b)\not=0$ and then
\begin{center}
$\AA^{-1}=\dfrac{1}{4b}\Circ\big(\b(2,1,0)\big)-\dfrac{a}{b(an+b)}\,\JJ.$
\end{center}
\end{itemize}
\end{corollary}
 \vspace{.25cm}

We end up this paper  by deriving the inverse of  a general symmetric circulant tridiagonal matrix,   without assuming the hypothesis of diagonally dominance.  Notice the difference between  our result and the methodology given in \cite{Ro90}.
%\cite{R90}.

\begin{corollary}
\label{rojo}
For $a,b\in \RR$, $b\not=0$, the circulant matrix $\Circ(a,b,0,\ldots,0,b)$ is invertible iff
$$\prod\limits_{j=0}^{\lceil\frac{n-1}{2}\rceil}\Big[a+2b\cos\Big(\frac{2\pi j}{n}\Big)\Big]\not=0$$
and, in this case
$$\Circ(a,b,0,\ldots,0,b)^{-1}=\Circ\big(\g(a,b,0)\big),$$ where
$$g_j(a,b,0)=\dfrac{(-1)^j}{2b[1-(-1)^nT_n(\frac{a}{2b})]}\Big[U_{j-2}\Big(\frac{a}{2b}\Big)+(-1)^{n}U_{n-j}\Big(\frac{a}{2b}\Big)\Big],\hspace{.25cm}j=1,\ldots,n.$$

\end{corollary}
Notice that the diagonally dominant hypothesis $|a|>2|b|$ clearly implies that $a+2b\cos\big(\frac{2\pi j}{n}\big)\not=0$ for any $j=0,\ldots,n$.

\vspace{.5cc}

\section*{Acknowledgments}
%{\bf Acknowledgments.}
This work has been partly supported by the Spanish Research Council (Comisi\'on Interministerial de Ciencia y Tecnolog\'\i a,)
under projects MTM2011-28800-C02-01 and MTM2011-28800-C02-02.

\vspace{.5cm}

\end{document}